\documentclass[12pt]{article}
\usepackage{amsmath,amsfonts,amsthm,mathrsfs}

\renewcommand{\epsilon}{\varepsilon } 
 
\renewcommand{\rho}{\varrho } 
\renewcommand{\phi}{\varphi }

\newtheorem{theorem}{Theorem}
\newtheorem{lemma}{Lemma}
\newtheorem{corollary}{Corollary}
\newtheorem{proposition}{Proposition}

\begin {document}
 \title{
Randomized Complexity of Vector-Valued Approximation
}

\author {Stefan Heinrich\\
Department of Computer Science\\
RPTU Kaiserslautern-Landau\\
D-67653 Kaiserslautern, Germany}  
\date{}
\maketitle

\begin{abstract} 
We study the  randomized $n$-th minimal errors (and hence the complexity) of  vector valued approximation.
In a recent paper by the author [Randomized complexity of parametric integration and the role of adaption I. Finite dimensional case (preprint)] a long-standing problem of Information-Based Complexity was solved: Is there a constant $c>0$ such that for all linear problems $\mathcal{P}$ the randomized non-adaptive and adaptive $n$-th minimal errors can deviate at most by a factor of $c$? That is, does the following hold for all linear $\mathcal{P}$ and $n\in {\mathbb N}$
\begin{equation*}
e_n^{\rm ran-non} (\mathcal{P})\le ce_n^{\rm ran} (\mathcal{P}) \, {\bf ?}
\end{equation*}
The analysis of vector-valued mean computation showed that the answer is negative. More precisely, there are instances of this problem where the gap between non-adaptive and adaptive randomized minimal errors can be (up to log factors) of the order $n^{1/8}$. This raises the question about the maximal possible deviation. In this paper we 
show that for certain instances of vector valued approximation the gap is $n^{1/2}$ (again, up to log factors).
\end{abstract}

\section{Introduction}
\label{sec:1}
Let $N,N_1,N_2\in{\mathbb{N}}$  and  $1\leq p,q,u,v \leq \infty$. We define the space
$L_p^N$ as the set of all functions $f: {\mathbb{Z}}[1,N]:=\{1,2,\dots,N\} \rightarrow {\mathbb{K}}$
with the norm
\begin{eqnarray*}
\| f \|_{L_p^N} = 
  \left( \frac{1}{N} \sum_{i=1}^N |f(i)|^p \right)^{1/p} \;(p<\infty),\quad  
  \| f \|_{L_\infty^N}=\max_{1\le i\le N} |f(i)| .
\end{eqnarray*}
and the space $L_p^{N_1}\big(L_u^{N_2}\big)$ as the set of all functions $f: {\mathbb{Z}}[1,N_1]\times {\mathbb{Z}}[1,N_2]\rightarrow {\mathbb{K}}$
with the norm
\[
\| f \|_{L_p^{N_1}\big(L_u^{N_2}\big)} = \Big\|\big(\|f_i\|_{L_u^{N_2}}\big)_{i=1}^{N_1}\Big\|_{L_p^{N_1}}
\]
with 
$
f_i=(f(i,j))_{j=1}^{N_2}
$
being the rows of the matrix $(f(i,j))$.
In the present paper we study the complexity of approximation in the randomized setting.
More precisely, we determine the order of the randomized $n$-th minimal errors of 
\begin{equation}
\label{C6}
J^{N_1,N_2}:L_p^{N_1}\big(L_u^{N_2}\big)\to L_q^{N_1}\big(L_v^{N_2}\big), \quad J^{N_1,N_2}f=f.
\end{equation}
The input set is the unit ball of $L_p^{N_1}\big(L_u^{N_2}\big)$, the error is measured in the norm of $L_q^{N_1}\big(L_v^{N_2}\big)$ and information is standard (values of $f$).

It is well-known since the 80ies that for linear problems adaptive and non-adaptive $n$-th minimal errors can deviate at most by a factor of 2, thus for any linear problem $\mathcal{P}=(F,G,S,K,\Lambda)$ and any $n\in {\mathbb{N}}$
\begin{equation}
\label{J9}
e_n^{\rm det-non } (S,F,G)\le 2e_n^{\rm det } (S,F,G),
\end{equation}
see Gal and Micchelli \cite{GM80}, Traub and Wo\'zniakowski \cite{TW80}. 
The randomized analogue of this problem is as follows: Is there a constant $c>0$ such that for all linear problems
$\mathcal{P}=(F,G,S,K,\Lambda)$ and all $n\in{\mathbb{N}}$
\begin{equation*}
e_n^{\rm ran-non} (S,F,G)\le ce_n^{\rm ran } (S,F,G) \, {\bf ?}
\end{equation*}
See the open problem on p.\ 213 of \cite{Nov96}, and  Problem 20 on p.\ 146 of \cite{NW08}. This problem was solved recently by the author in \cite{Hei23a}, where it was shown that for some instances of vector-valued mean computation the gap between non-adaptive and adaptive randomized $n$-th minimal errors can be (up to log factors) of  order $n^{1/8}$. This raises the question about the maximal possible deviation. In this paper we study the randomized complexity of vector valued approximation and
show that for certain instances the gap is $n^{1/2}$ (again, up to log factors), see Corollary \ref{cor:1}.

\section{Preliminaries}
\label{sec:2}
Throughout this paper $\log$ means $\log_2$. We denote ${\mathbb{N}}=\{1,2,\dots\}$ and ${\mathbb{N}}_0={\mathbb{N}}\cup\{0\}$. The symbol ${\mathbb{K}}$ stands for the scalar field, which is either ${\mathbb R}$ or ${\mathbb C}$.
We often use the same symbol
$c, c_1,c_2,\dots$ for possibly different constants, even if they appear in a sequence
of relations. However, some constants are supposed to have the same meaning throughout a proof  --  these are denoted by symbols $c(1),c(2),\dots$. The unit ball of a normed space $X$ is denoted by $B_X$. 

We work in the framework of IBC \cite{Nov88,TWW88}, using specifically the general approach from \cite{Hei05a, Hei05b}, see also the extended introduction in \cite{Hei23a}. We refer to these papers for notation and background.

An abstract  numerical problem $\mathcal{P}$ is given as 
$\mathcal{P}=(F,G,S,K,\Lambda)$, where
 $F$ is a non-empty set, 
$G$ a Banach space, and $S$ is a mapping $F\to G$.  The operator $S$ is called the  solution operator, it sends the input  $f\in F$ of our problem to the exact solution $S(f)$. Moreover, $\Lambda$ is a nonempty set of mappings from $F$ to $K$, the set of information functionals, where $K$ is any nonempty set -- the set of values of information functionals. 

A problem $\mathcal{P}$ is called linear, if $K={\mathbb K}$, $F$ is a convex and balanced subset of a linear space $X$ over ${\mathbb K}$,  
$S$ is the restriction to $F$ of a linear operator 
from $X$ to $G$, and each $\lambda\in\Lambda$ is the restriction  to $F$ of a linear mapping from $X$ to ${\mathbb K}$.

In this paper we consider the linear problem
$$
\mathcal{P}^{N_1,N_2}=\left( B_{L_p^{N_1}\big(L_u^{N_2}\big)},L_q^{N_1}(L_v^{N_2}),J^{N_1,N_2},{\mathbb{K}},\Lambda\right),
$$
where
$\Lambda=\{\delta_{ij}:\, 1\le i\le N_1,\,1\le j\le N_2\}$ with $\delta_{ij}(f)=f(i,j)$. 

A deterministic algorithm for $\mathcal{P}$ 
is a tuple $A=((L_i)_{i=1}^\infty, (\tau_i)_{i=0}^\infty,(\varphi_i)_{i=0}^\infty)$
such that 
$L_1\in\Lambda$, $\tau_0\in\{0,1\}$, $\varphi_0\in G$,
and for $i\in {\mathbb{N}}$,
$L_{i+1} : K^i\to \Lambda$,
$\tau_i:  K^i\to \{0,1\}$, and 
$\varphi_i:  K^i\to G$ 
are arbitrary mappings, where $K^i$ denotes the $i$-th Cartesian power of $K$.
Given an input $f\in F$, we define $(\lambda_i)_{i=1}^\infty$ with $\lambda_i\in \Lambda$ 
as follows:
\begin{eqnarray*}
\lambda_1=L_1, \quad
\lambda_i=L_i(\lambda_1(f),\dots,\lambda_{i-1}(f))\quad(i\ge 2)\label{RC2}.
\end{eqnarray*}
Define ${\rm card}(A,f)$, the cardinality  of $A$ at input $f$, to be $0$ if $\tau_0=1$. If $\tau_0=0$, let ${\rm card}(A,f)$ be
the first integer $n\ge 1$ with $
\tau_n(\lambda_1(f),\dots,\lambda_n(f))=1 $
if there is such an $n$. If $\tau_0=0$ and no such $n\in {\mathbb{N}}$ exists, 
put ${\rm card}(A,f)=+\infty$.
We define the output $A(f)$ of algorithm $A$ at input $f$ as
\begin{equation}
\label{K8}
A(f)=\left\{\begin{array}{lll}
\varphi_0  & \mbox{if} \quad {\rm card}(A,f)\in \{0,\infty\} \\[.2cm]
\varphi_n(\lambda_1(f),\dots,\lambda_n(f))  &\mbox{if} \quad 1\le {\rm card}(A,f)=n<\infty. 
\end{array}
\right.
\end{equation}
The cardinality of $A$ is defined as
$
{\rm card}(A,F)=\sup_{f\in F}{\rm card} (A,f).
$
Given $n\in{\mathbb{N}}_0$, we define $\mathscr{A}_n^{\rm det }(\mathcal{P})$ as the set of 
deterministic algorithms $A$ for $\mathcal{P}$ with ${\rm card} (A)\le n$
and the deterministic $n$-th minimal error of $S$ as
\begin{equation}
\label{X1}
e_n^{\rm det } (S,F,G)=\inf_{A\in\mathscr{A}_n^{\rm det }(\mathcal{P}) }  
\sup_{f\in F}\|S(f)-A(f)\|_G.
\end{equation}
A deterministic algorithm is called non-adaptive, if all $L_i$ and all $\tau_i$ are constant, in other words, 
$L_i\in \Lambda$, $\tau_i\in\{0,1\}$.
The subset of non-adaptive algorithms in  $\mathscr{A}_n^{\rm det }(\mathcal{P})$ is denoted by  $\mathscr{A}_n^{\rm det-non } (\mathcal{P})$ and the non-adaptive deterministic $n$-th minimal error $e_n^{\rm det-non} (S,F,G)$  is defined in analogy with \eqref{X1}.

A randomized algorithm for $\mathcal{P}$ is a tuple
$
A = (( \Omega, \Sigma, {\mathbb P}),(A_{\omega})_{\omega \in \Omega}),
$
where $(\Omega, \Sigma, {\mathbb P})$ is a probability space and for each $\omega \in \Omega$, $A_{\omega}$ is a deterministic algorithm for $\mathcal{P}$. 
Let $n\in{\mathbb{N}}_0$. Then $\mathscr{A}_n^{{\rm ran }}(\mathcal{P})$ stands for the class of 
randomized algorithms $A$ for $\mathcal{P}$ with the following properties: For each $f\in F$ the mapping $\omega\to {\rm card}(A_{\omega},f)$ is $\Sigma$-measurable,  
$
{\mathbb E}\,{\rm card} (A_{\omega},f)\le n,
$
and the mapping $\omega\to A_{\omega}(f)$ is $\Sigma$-to-Borel measurable  and ${\mathbb P}$-almost surely separably valued, i.e.,  there is a separable subspace $G_f$ of $G$ such that 
${\mathbb P}\{\omega:\,A_\omega(f)\in G_f\}=1$.
We define the  cardinality of $A\in \mathscr{A}_n^{{\rm ran }}(\mathcal{P})$ as
$
{\rm card}(A,F)=\sup_{f\in F}{\mathbb E}\,{\rm card} (A_{\omega},f),
$
and the randomized $n$-th minimal error of $S$ as
\begin{equation*}
e_n^{\rm ran } (S,F,G)=\inf_{A\in\mathscr{A}_n^{\rm ran }(\mathcal{P}) }  
\sup_{f\in F}{\mathbb E}\,\|S(f)-A_{\omega}(f)\|_G.
\end{equation*}
We call a randomized algorithm $(( \Omega, \Sigma, {\mathbb P}), (A_{\omega})_{\omega \in \Omega})$ non-adaptive, if $A_{\omega}$ is non-adaptive for all $\omega \in \Omega$. Furthermore, $\mathscr{A}_n^{\rm ran-non}(\mathcal{P})$ is the subset of $\mathscr{A}_n^{\rm ran }(\mathcal{P})$ consisting of non-adaptive algorithms,
and $e_n^{\rm ran-non} (S,F,G)$ denotes the non-adaptive randomized $n$-th minimal error.

We also need the average case setting. For the purposes of this paper we consider it only for measures which are supported by a finite subset of $F$. Then the underlying $\sigma$-algebra is assumed to be $2^F$, therefore no measurability conditions have to be imposed on $S$ and the involved   deterministic algorithms. So let $\mu$ be a  probability measure on $F$ with finite support, let ${\rm card} (A, \mu)  =  \int_{F} {\rm card}(A,f )  d \mu(f)$,
and define 
\begin{eqnarray*}
e_n^{\rm avg } (S, \mu,G) &=& \inf_A\int_{F} \| S(f) 
	- A(f) \|_G  d \mu(f),
\end{eqnarray*}
where the infimum is taken over all deterministic algorithms with ${\rm card} (A, \mu)\le n$. Correspondingly, $e_n^{\rm avg-non} (S, \mu,G)$ is defined. 
We use the following well-known results to prove lower bounds. 
\begin{lemma}
\label{Ulem:5}
For every probability measure $\mu$ on $F$ of finite support we have 
\begin{eqnarray*}
e_n^{\rm ran }(S,F)\ge \frac{1}{2}e_{2n}^{\rm avg }(S,\mu),\quad
e_n^{\rm ran-non}(S,F)\ge \frac{1}{2}e_{2n}^{\rm avg-non}(S,\mu).
\end{eqnarray*}
\end{lemma}
The types of lower bounds stated in the next lemma are well-known in IBC (see \cite{Nov88,TWW88}). 
For the specific form  presented here we refer, e.g., to \cite{Hei05a},  Lemma 6 for statement (i), and to \cite{Hei18a}, Proposition 3.1 for (ii). 
\begin{lemma}\label{lem:5} 
Let $\mathcal{P}=(F,G,S,K,\Lambda)$ be a linear problem, 
$\bar{n}\in{\mathbb{N}}$, and suppose there are $(f_i)_{i=1}^{\bar{n}}\subseteq F$
such that the sets $\{\lambda\in \Lambda\,:\, \lambda(f_i)\ne 0\}\; (i=1,\dots,\bar{n})$
are mutually disjoint.
Then the following hold for all $n\in{\mathbb{N}}$ with
$4n<\bar{n}$:\\
\indent (i) If $\sum_{i=1}^{\bar{n}} \alpha_i f_i\in F$ for all sequences $(\alpha_i)_{i=1}^{\bar{n}}\in \{-1,1\}^{\bar{n}}$ and $\mu$ is the distribution of $\sum_{i=1}^{\bar{n}} \varepsilon_i f_i$, where $\varepsilon_i$ are independent Bernoulli random variables
with ${\mathbb P}\{\varepsilon_i=1\}={\mathbb P}\{\varepsilon_i=-1\}=1/2$,  then   
$$
e_n^{\rm avg }(S,\mu)\ge \frac{1}{2}\min\bigg\{{\mathbb E}\,\Big\|\sum_{i\in I}\varepsilon_iSf_i\Big\|_G:\,I\subseteq\{1,\dots,\bar{n}\},\,|I|\ge \bar{n}-2n\bigg\}.
$$
\indent (ii) If $\alpha f_i\in F$ for all $1\le i\le\bar{n}$ and $\alpha\in \{-1,1\}$, and $\mu$ is the uniform distribution on the set $\{\alpha f_i\,:\,1\le i\le\bar{n},\; \alpha\in \{-1,1\}\}$,
then 
$$
e_n^{\rm avg }(S,\mu)\ge \frac{1}{2}\min_{1\le i\le \bar{n}} \|Sf_i\|_G.
$$
\end{lemma}
Finally, let $\theta$ be the mapping given by the median, that is, if $z_1^*\le  \dots\le z_m^*$ is the non-decreasing rearrangement of $(z_1,\dots,z_m)\in{\mathbb R}^m$, then  
$\theta(z_1,\dots,z_m)$ stands for 
$z^*_{(m+1)/2}$ if $m$ is odd and $\frac{z_{m/2}^*+z_{m/2+1}^*}{2}$ if $m$ is even. The following is well-known,  see, e.g, \cite{{Hei01}}.
\begin{lemma}
\label{Ulem:2e}
Let $\zeta_1,\dots,\zeta_m$ be independent, identically distributed real-valued random variables on a probability space $(\Omega,\Sigma,{\mathbb P})$, $z\in {\mathbb R}$, $\varepsilon>0$ , and assume that ${\mathbb P}\{|z-\zeta_1|\le\varepsilon\}\ge 3/4$. Then
\begin{equation*}
{\mathbb P}\{|z-\theta(\zeta_1,\dots,\zeta_m)|\le \varepsilon\}\ge 1-e^{-m/8}.
\end{equation*}
\end{lemma}

\section{An adaptive algorithm for vector valued approximation}
\label{sec:4}
We refer to the definition of the embedding $J^{N_1,N_2}$ given in \eqref{C6}.
It is easily checked by  H\"older's inequality that 
\begin{equation}
\label{WN2}
\big\|J^{N_1,N_2}\big\|=N_1^{\left(1/p-1/q\right)_+}N_2^{\left(1/u-1/v\right)_+},
\end{equation}
with $a_+=\max(a,0)$ for $a\in {\mathbb R}$. We need the randomized norm estimation algorithm from \cite{Hei18}. 
 Let $(Q,\mathcal{Q},\varrho)$ be a probability space and let $1\le v< u
\le \infty$. 
For $n\in{\mathbb{N}}$ define  $A_n^1=(A_{n,\omega}^1)_{\omega\in\Omega}$ by setting 
 for $\omega\in \Omega$ and $f\in L_u(Q,\mathcal{Q},\varrho)$ 
\begin{eqnarray}
A_{n,\omega}^1(f)&=&\left(\frac{1}{n}\sum_{i=1}^n|f(\xi_i(\omega_2))|^v\right)^{1/v},\label{QG4c}
\end{eqnarray}
where $\xi_i$ are independent $Q$-valued random variables on a probability space $(\Omega,\Sigma,{\mathbb P})$ with distribution $\varrho$. 
The following is essentially  Proposition 6.3 of \cite{Hei18}, for a self-contained proof we refer to \cite{Hei23a}. 
\begin{proposition}
\label{Rpro:4}
Let $1\le v<u\le \infty$.  
Then there is 
a constant $c>0$ such that for  all probability spaces $(Q,\mathcal{Q},\varrho)$, $f\in L_u(Q,\mathcal{Q},\varrho)$, and  $n\in{\mathbb{N}}$ 
\begin{eqnarray}
{\mathbb E}\,\left|\|f\|_{L_v(Q,\mathcal{Q},\varrho)}-A_{n,\omega}^1(f)\right| &\le& cn^{\max\left(1/u-1/v,-1/2\right)}\|f\|_{L_u(Q,\mathcal{Q},\varrho)}.
\label{RH4}
\end{eqnarray}
\end{proposition}
The algorithm for approximation of $J^{N_1,N_2}$ will only be defined for the case that $1\le p<q\le \infty$ and $1\le v<u\le \infty$ (it turns out that for the other cases the zero algorithm is of optimal order). Define for  $m,n\in{\mathbb{N}}$, $n<N_1N_2$  an adaptive algorithm $A_{n,m,\omega}^2$. Let $f\in L_p^{N_1}(L_u^{N_2})$ and set $f_i=(f(i,j))_{j=1}^{N_2}$.  Let
$
\left\{\xi_{jk}:\, 1\le j\le \left\lceil\frac{n}{N_1}\right\rceil,\,1\le k\le m\right\}
$
be independent random variables on a probability space $(\Omega,\Sigma,{\mathbb P})$ uniformly distributed over $\{1,\dots,N_2\}$. 

We apply algorithm $A_{n,\omega}^1$, see \eqref{QG4c}, to estimate $\|f_i\|_{L_v^{N_2}}$ by setting  for $\omega\in\Omega$, $1\le i\le N_1$, $1\le k\le m$  
\begin{eqnarray*}
a_{ik}(\omega)&=&\Bigg(\left\lceil\frac{n}{N_1}\right\rceil^{-1}\sum_{1\le j\le \left\lceil\frac{n}{N_1}\right\rceil} |f_i(\xi_{jk}(\omega))|^v\Bigg)^{1/v}, \quad
\tilde{a}_i(\omega)=\theta\big((a_{ik}(\omega))_{k=1}^{m}\big).
\end{eqnarray*}
Let $\tilde{a}_{\pi(1)}\ge \dots\ge \tilde{a}_{\pi(N_1)}$ be a non-increasing rearrangement of $(\tilde{a}_i)$, with $\pi$ a permutation.
Then the output $A_{n,m,\omega}^2(f)\in L_q^{N_1}(L_v^{N_2}) $ of the algorithm is defined as 
\begin{equation}
\label{UN3}
A_{n,m,\omega}^2(f)=b=(b_i(\omega))_{i=1}^{N_1},\quad b_{\pi(i)}(\omega)=\left\{\begin{array}{lll}
  f_{\pi(i)} & \quad\mbox{if}\quad i\le \left\lceil\frac{n}{N_2}\right\rceil   \\
 0 &  \quad\mbox{otherwise}
    \end{array}
\right. 
\end{equation}
(note that the assumption on $n$ implies  $ \left\lceil\frac{n}{N_2}\right\rceil\le N_1$).
 If $\frac{1}{v}-\frac{1}{u}> \frac{1}{2}$, we use an iterated  version $(A_{n,m,\omega}^3)_{\omega\in \Omega}$, where $(\Omega,\Sigma,{\mathbb P})=(\Omega_1,\Sigma_1,{\mathbb P}_1)\times (\Omega_2,\Sigma_2,{\mathbb P}_2)$ with $(\Omega_\iota,\Sigma_\iota,{\mathbb P}_\iota)$ $(\iota=1,2)$ being probability spaces.
We define
\begin{equation}
\label{H3}
A_{n,m,\omega}^3(f)=A_{n,m,\omega_1}^2(f)+A_{n,m,\omega_2}^2(f-A_{n,m,\omega_1}^2(f))\quad (\omega=(\omega_1,\omega_2)).
\end{equation}
\indent The constants in the subsequent statements and proofs are independent of the parameters $n$, $N_1$,$N_2$, and $m$. This is also made clear by the order of quantifiers in the respective statements.

\begin{proposition}
\label{pro:5}
Let $1\le p<q\le \infty$, $1\le v<u\le \infty$, and  $1\le w<\infty$. Then there exist constants $c_1>1$, $c_2>0$ such that the following hold for all $m,n,N_1,N_2\in{\mathbb{N}}$ with  $n<N_1N_2$ and $ f\in L_p^{N_1}(L_u^{N_2})$:
\begin{equation}
\label{WM1}
{\rm card}(A_{n,m,\omega}^2)\le (m+1)n+mN_1+N_2, \quad {\rm card}(A_{n,m,\omega}^3)=2\,{\rm card}(A_{n,m,\omega}^2).
\end{equation}
Furthermore, setting $A_{n,m,\omega}=A_{n,m,\omega}^2$ if $\frac{1}{v}-\frac{1}{u}\le \frac{1}{2} $ and $A_{n,m,\omega}=A_{n,m,\omega}^3$ if $\frac{1}{v}-\frac{1}{u}> \frac{1}{2} $,
we have for $m\ge c_1\log(N_1+N_2)$
\begin{eqnarray}
\label{mc-eq:2}
\lefteqn{\left({\mathbb E}\, \|f-A_{n,m,\omega}(f)\|_{L_q^{N_1}(L_v^{N_2})}^{w}\right)^{1/w}}
\nonumber\\
&\le& 
 c_2N_1^{1/p-1/q}\Bigg(\left\lceil\frac{n}{N_1}\right\rceil^{1/u-1/v}+ \left\lceil\frac{n}{N_2}\right\rceil^{1/q-1/p}\Bigg)\|f\|_{L_p^{N_1}(L_u^{N_2})}.
\end{eqnarray}
\end{proposition}
\begin{proof}
The total number of samples in $A_{n,m,\omega}^2$ is 
\begin{eqnarray*}
mN_1\left\lceil\frac{n}{N_1}\right\rceil+N_2\left\lceil\frac{n}{N_2}\right\rceil
&\le& mn+mN_1+n+N_2,
\end{eqnarray*}
which gives \eqref{WM1}. For $n<\max(N_1,N_2)$ relation \eqref{mc-eq:2} follows from \eqref{WN2}. Hence in the sequel we assume $n\ge \max(N_1,N_2)$. 
Fix $f\in L_p^{N_1}(L_u^{N_2})$. First we consider the case $\frac{1}{v}-\frac{1}{u}\le \frac{1}{2}$. By Proposition \ref{Rpro:4}, where here the respective constant is denoted by $c(0)$,
\begin{equation}
\label{C0}
{\mathbb E}\,\Big|\|f_i\|_{L_v^{N_2}}-a_{ik}\Big|\le c(0) \bigg(\frac{n}{N_1}\bigg)^{1/u-1/v}\|f_i\|_{L_u^{N_2}},
\end{equation}
and therefore, 
\begin{equation}
\label{A6}
{\mathbb P}\left\{\omega\in \Omega:\,\Big|\|f_i\|_{L_v^{N_2}}-a_{ik}(\omega)\Big|\le 4c(0) \bigg(\frac{n}{N_1}\bigg)^{1/u-1/v}\|f_i\|_{L_u^{N_2}}\right\}\ge \frac{3}{4}.
\end{equation}
Let 
$
c(1)=\frac{8(w+1)}{\log e}>1
$
(recall that $\log$ always means $\log_2$), then $m\ge c(1)\log(N_1+N_2)$ implies $e^{-m/8}\le (N_1+N_2)^{-w-1}$.
From \eqref{A6} and Lemma \ref{Ulem:2e} we conclude
\begin{equation*}
{\mathbb P}\left\{\omega\in \Omega:\,\Big|\|f_i\|_{L_v^{N_2}}-\tilde{a}_i(\omega)\Big|\le 4c(0) \bigg(\frac{n}{N_1}\bigg)^{1/u-1/v}\|f_i\|_{L_u^{N_2}}\right\}\ge 1-(N_1+N_2)^{-w-1}.
\end{equation*}
Let
\begin{eqnarray}
\Omega_0&=&\Bigg\{\omega\in \Omega:
\Big|\|f_i\|_{L_v^{N_2}}-\tilde{a}_i(\omega)\Big|\le 4c(0) \bigg(\frac{n}{N_1}\bigg)^{1/u-1/v}\|f_i\|_{L_u^{N_2}}\;(1\le i\le N_1)\Bigg\},\label{L0}
\end{eqnarray}
thus
\begin{equation}
\label{L1}
{\mathbb P}(\Omega_0)\ge 1-(N_1+N_2)^{-w}.
\end{equation}
 Fix $\omega\in \Omega_0$.
Then by \eqref{L0} for all $i$
\begin{equation}
\label{UL7}
\tilde{a}_i(\omega)\le c \|f_i\|_{L_u^{N_2}}.
\end{equation}
Consequently,
\begin{equation}
\label{UL9}
\bigg(\frac{1}{N_1}\sum_{i=1}^{N_1}\tilde{a}_i(\omega)^p\bigg)^{1/p}
\le c\bigg(\frac{1}{N_1}\sum_{i=1}^{N_1}\|f_i\|_{L_u^{N_2}}^p\bigg)^{1/p}= c\|f\|_{L_p^{N_1}(L_u^{N_2})}.
\end{equation}
Let $M=\left\lceil\frac{n}{N_2}\right\rceil$. It follows that
\begin{eqnarray}
c\|f\|_{L_p^{N_1}(L_u^{N_2})}\ge \bigg(\frac{1}{N_1}\sum_{i=1}^M\tilde{a}_{\pi(i)}^p\bigg)^{1/p}\ge \left(\frac{M}{N_1}\right)^{1/p}\tilde{a}_{\pi(M)},
\label{UM1}
\end{eqnarray}
thus for $i> M$
\begin{eqnarray}
\tilde{a}_{\pi(i)}\le \tilde{a}_{\pi(M)}\le c\left(\frac{N_1}{M}\right)^{1/p}\|f\|_{L_p^{N_1}(L_u^{N_2})}\le c\left(\frac{N_1N_2}{n}\right)^{1/p}\|f\|_{L_p^{N_1}(L_u^{N_2})}
\label{C8}.
\end{eqnarray}
Furthermore, by \eqref{UN3}, for $i\le M$
\begin{eqnarray}
\label{U8}
\big\|f_{\pi(i)}-b_{\pi(i)}(\omega)\big\|_{L_v^{N_2}}    &=& 0.
\end{eqnarray}
Let 
\begin{equation}
\label{UM3}
I(\omega):=\bigg\{1\le i\le N_1:\, \tilde{a}_{\pi(i)}(\omega)\le \frac{\|f_{\pi(i)}\|_{L_v^{N_2}}}{2}\bigg\},
\end{equation}
hence we conclude from  \eqref{UN3} and \eqref{L0} for $i>M,\,i\in I(\omega)$
\begin{eqnarray}
\big\|f_{\pi(i)}-b_{\pi(i)}(\omega)\big\|_{L_v^{N_2}}&=& \|f_{\pi(i)}\|_{L_v^{N_2}}
\le 2\big(\|f_{\pi(i)}\|_{L_v^{N_2}}-\tilde{a}_{\pi(i)}(\omega)\big) \nonumber\\
&\le&c\bigg(\frac{n}{N_1}\bigg)^{1/u-1/v}\|f_{\pi(i)}\|_{L_u^{N_2}}. \label{UL8}
\end{eqnarray}
On the other hand,  we have by \eqref{C8} for $i>M,\,i\not\in I(\omega)$
\begin{eqnarray}
\big\|f_{\pi(i)}-b_{\pi(i)}(\omega)\big\|_{L_v^{N_2}}& =&\big\|f_{\pi(i)}\big\|_{L_v^{N_2}} < 2\tilde{a}_{\pi(i)}(\omega)=
 2\tilde{a}_{\pi(i)}(\omega)^{p/q}\tilde{a}_{\pi(i)}(\omega)^{1-p/q}
\nonumber\\
&\le& c\tilde{a}_{\pi(i)}(\omega)^{p/q}  \left(\frac{N_1N_2}{n}\right)^{1/p-1/q}\|f\|_{L_p^{N_1}(L_u^{N_2})}^{1-p/q} 
\label{UM4} 
\end{eqnarray}
(with the convention $0^0=1$). Combining \eqref{U8}, \eqref{UL8}, and \eqref{UM4}, we get for $1\le i\le N_1$
\begin{eqnarray*}
\big\|f_i-b_i(\omega)\big\|_{L_v^{N_2}}
&\le& 
  c\bigg(\frac{n}{N_1}\bigg)^{1/u-1/v}\|f_i\|_{L_u^{N_2}}+c\tilde{a}_i(\omega)^{p/q}  \left(\frac{N_1N_2}{n}\right)^{1/p-1/q}\|f\|_{L_p^{N_1}(L_u^{N_2})}^{1-p/q}.
\end{eqnarray*}
Together with \eqref{UL9} we obtain for $\omega\in \Omega_0$, 
\begin{eqnarray}
\lefteqn{\big\|f-(b_i(\omega))_{i=1}^{N_1}\big\|_{L_q^{N_1}(L_v^{N_2})}}
\notag\\
&\le& 
c\bigg(\frac{n}{N_1}\bigg)^{1/u-1/v}\Big\|\big(\|f_i\|_{L_u^{N_2}}\big)_{i=1}^{N_1}\Big\|_{L_q^{N_1}}
+ c \left(\frac{N_1N_2}{n}\right)^{1/p-1/q}\Big\|\big(\tilde{a}_i(\omega)^{p/q}\big)_{i=1}^{N_1}\Big\|_{L_q^{N_1}} \|f\|_{L_p^{N_1}(L_u^{N_2})}^{1-p/q}
\notag\\
&\le& 
c\bigg(\frac{n}{N_1}\bigg)^{1/u-1/v}N_1^{1/p-1/q}\Big\|\big(\|f_i\|_{L_u^{N_2}}\big)_{i=1}^{N_1}\Big\|_{L_p^{N_1}}
\nonumber\\
&& 
+ c \left(\frac{N_1N_2}{n}\right)^{1/p-1/q}\Bigg(\frac{1}{N_1}\sum_{i=1}^{N_1}\tilde{a}_i(\omega)^p \Bigg)^{1/q}\|f\|_{L_p^{N_1}(L_u^{N_2})}^{1-p/q}
\nonumber\\
&\le&cN_1^{1/p-1/q}\Bigg(\bigg(\frac{n}{N_1}\bigg)^{1/u-1/v}+ \left(\frac{n}{N_2}\right)^{1/q-1/p}\Bigg)\|f\|_{L_p^{N_1}(L_u^{N_2})}.
\label{UN1}
\end{eqnarray}

To estimate the error on $\Omega\setminus \Omega_0$ we note that by \eqref{UN3} for all $\omega\in\Omega$, $b_i$ is either $f_i$ or zero.
Consequently
\begin{equation}
\label{M0}
\big\|f-b(\omega)\big\|_{L_q^{N_1}(L_v^{N_2})}\le\big\|f\big\|_{L_q^{N_1}(L_v^{N_2})}\le N_1^{1/p-1/q} \|f\|_{L_p^{N_1}(L_u^{N_2})},
\end{equation}
and therefore, using \eqref{L1},
\begin{eqnarray*}
\lefteqn{\left(\int_{\Omega\setminus \Omega_0}\big\|f-(b_i(\omega))_{i=1}^{N_1}\big\|_{L_q^{N_1}(L_v^{N_2})}^{w}d{\mathbb P}(\omega)\right)^{1/w}}
\\
&\le& N_1^{1/p-1/q}(N_1+N_2)^{-1}\|f\|_{L_p^{N_1}(L_u^{N_2})}
\le N_1^{1/p-1/q}\bigg(\frac{n}{N_1}\bigg)^{1/q-1/p}\|f\|_{L_p^{N_1}(L_u^{N_2})},
\end{eqnarray*}
the last relation being a consquence of $n< N_1N_2$. Together with \eqref{UN1} this shows \eqref{mc-eq:2} under the assumption $\frac{1}{v}-\frac{1}{u}\le \frac{1}{2}$.

Finally we consider the case $\frac{1}{v}-\frac{1}{u}>\frac{1}{2}$. We define $q_1, v_1$ by
\begin{equation}
\label{H0}
\frac{1}{q_1}=\frac{1}{2}\left(\frac{1}{p}+\frac{1}{q}\right), \quad \frac{1}{v_1}=\frac{1}{2}\left(\frac{1}{u}+\frac{1}{v}\right),
\end{equation}
then $1\le p<q_1<q$, $v<v_1<u$,  $\frac{1}{v_1}-\frac{1}{u}\le \frac{1}{2}$, and $\frac{1}{v}-\frac{1}{v_1}\le\frac{1}{2}$, so we conclude from the already shown case of
\eqref{mc-eq:2}
\begin{eqnarray}
\label{H1}
\lefteqn{\left({\mathbb E}\,_{\omega_1} \|f-A_{n,m,\omega_1}^2(f)\|_{L_{q_1}^{N_1}(L_{v_1}^{N_2})}^{w}\right)^{1/w}}
\nonumber\\
&\le& 
 cN_1^{1/p-1/q_1}\Bigg(\bigg(\frac{n}{N_1}\bigg)^{1/u-1/v_1}+ \left(\frac{n}{N_2}\right)^{1/q_1-1/p}\Bigg)\|f\|_{L_p^{N_1}(L_u^{N_2})}
\end{eqnarray}
and, with $g= f-A_{n,m,\omega_1}^2(f)$,
\begin{eqnarray}
\label{H2}
\lefteqn{\left({\mathbb E}\,_{\omega_2} \|g-A_{n,m,\omega_2}^2(g)\|_{L_q^{N_1}(L_v^{N_2})}^{w}\right)^{1/w}}
\nonumber\\
&\le& 
cN_1^{1/q_1-1/q}\Bigg(\bigg(\frac{n}{N_1}\bigg)^{1/v_1-1/v}+ \left(\frac{n}{N_2}\right)^{1/q-1/q_1}\Bigg)\|g\|_{L_{q_1}^{N_1}(L_{v_1}^{N_2})}.
\end{eqnarray}
From  \eqref{H3}, \eqref{H1}, and \eqref{H2} we obtain 
%
\begin{eqnarray}
\label{H4}
\lefteqn{\left({\mathbb E}\,_{\omega_1}{\mathbb E}\,_{\omega_2} \|f-A_{n,m,(\omega_1,\omega_2)}^3(f)\|_{L_q^{N_1}(L_v^{N_2})}^{w}\right)^{1/w}}
\nonumber\\
&=& \left({\mathbb E}\,_{\omega_1}{\mathbb E}\,_{\omega_2} \|f-A_{n,m,\omega_1}^2(f)-A_{n,m,\omega_2}^2(f-A_{n,m,\omega_1}^2(f))\|_{L_q^{N_1}(L_v^{N_2})}^{w}\right)^{1/w}
\nonumber\\
&\le& 
cN_1^{1/q_1-1/q}\Bigg(\bigg(\frac{n}{N_1}\bigg)^{1/v_1-1/v}+ \left(\frac{n}{N_2}\right)^{1/q-1/q_1}\Bigg)\left({\mathbb E}\,_{\omega_1}\|f-A_{n,m,\omega_1}^2(f)\|_{L_{q_1}^{N_1}(L_{v_1}^{N_2})}^w\right)^{1/w}
\nonumber
\\
&\le& 
cN_1^{1/q_1-1/q}\Bigg(\bigg(\frac{n}{N_1}\bigg)^{1/v_1-1/v}+ \left(\frac{n}{N_2}\right)^{1/q-1/q_1}\Bigg)
\nonumber\\
&&\times N_1^{1/p-1/q_1}\Bigg(\bigg(\frac{n}{N_1}\bigg)^{1/u-1/v_1}+ \left(\frac{n}{N_2}\right)^{1/q_1-1/p}\Bigg)\|f\|_{L_p^{N_1}(L_u^{N_2})}
\nonumber
\\
&=&cN_1^{1/p-1/q}\Bigg(\bigg(\frac{n}{N_1}\bigg)^{\frac{1}{2}\left(1/u-1/v\right)}+ \left(\frac{n}{N_2}\right)^{\frac{1}{2}\left(1/q-1/p\right)}\Bigg)^2\|f\|_{L_p^{N_1}(L_u^{N_2})}
\nonumber\\
&\le&cN_1^{1/p-1/q}\Bigg(\bigg(\frac{n}{N_1}\bigg)^{1/u-1/v}+ \left(\frac{n}{N_2}\right)^{1/q-1/p}\Bigg)\|f\|_{L_p^{N_1}(L_u^{N_2})}.
\end{eqnarray}
This gives \eqref{mc-eq:2} and concludes the proof.
\end{proof}
\section{Lower bounds and complexity}
\begin{proposition}
\label{pro:1}
Let $1\le p,q,u,v\le \infty$. Then there exist constants $0<c_0<1$, $c_1\dots c_6>0$ such that for all $n,N_1,N_2\in {\mathbb{N}}$, with $n<c_0N_1N_2$ there exist probability measures $\mu_{n,N_1,N_2}^{(i)}$ $(1\le i\le 6)$ with finite support in $B_{L_p^{N_1}(L_u^{N_2})}$ such that 
\begin{eqnarray}
e_{n}^{\rm avg }(J^{N_1,N_2},\mu_{n,N_1,N_2}^{(1)},L_q^{N_1}(L_v^{N_2}))
&\ge&  c_1N_1^{1/p-1/q} \left\lceil\frac{n}{N_1}\right\rceil^{1/u-1/v}
\label{UL2}\\
e_{n}^{\rm avg }(J^{N_1,N_2},\mu_{n,N_1,N_2}^{(2)},L_q^{N_1}(L_v^{N_2}))
&\ge&  c_2N_1^{1/p-1/q} \left\lceil\frac{n}{N_2}\right\rceil^{1/q-1/p}
\label{L5}\\
e_n^{\rm avg }(J^{N_1,N_2},\mu_{n,N_1,N_2}^{(3)},L_q^{N_1}(L_v^{N_2}))&\ge&  c_4N_1^{1/p-1/q}N_2^{1/u-1/v}
\label{N3}\\
e_{n}^{\rm avg }(J^{N_1,N_2},\mu_{n,N_1,N_2}^{(4)},L_q^{N_1}(L_v^{N_2}))
&\ge& c_3 
\label{UL1}\\
e_{n}^{\rm avg }(J^{N_1,N_2},\mu_{n,N_1,N_2}^{(5)},L_q^{N_1}(L_v^{N_2}))
&\ge& c_5 \left\lceil\frac{n}{N_1}\right\rceil^{1/u-1/v}
\label{UL3}\\
e_{n}^{\rm avg-non}\big(J^{N_1,N_2},\mu_{n,N_1,N_2}^{(6)}\big),L_q^{N_1}(L_v^{N_2}))
&\ge& c_6 N_1^{1/p-1/q}.
\label{UQ3}
\end{eqnarray}
\end{proposition}
\begin{proof} We set $c_0=\frac1{21}$ and let $n\in {\mathbb{N}}$ be such that 
\begin{equation}
\label{B1}
1\le n< \frac{N_1N_2}{21}.
\end{equation}
Define for $L$ with  $1\le L\le N_2$ disjoint subsets of $\{1,\dots,N_2\}$ by setting
\begin{equation}
\label{UWL7}
D_{j}= \left\{ (j-1)\left\lfloor \frac{N_2}{L}\right\rfloor+1, \dots, 
j \left\lfloor \frac{N_2}{L}\right\rfloor  \right\} , \quad (j=1, \dots, L),
\end{equation}
then
\begin{equation}
\label{UL0}
 \frac{N_2}{2L}<  \left\lfloor \frac{N_2}{L}\right\rfloor=|D_{j}| \le \frac{N_2}{L}.
\end{equation}
To show \eqref{UL2}, we put 
\begin{equation}
\label{UK5}
L=\left\lfloor\frac{4n}{N_1}\right\rfloor +1,
\end{equation}
thus
\begin{equation}
\label{B6}
\left\lceil\frac{4n}{N_1}\right\rceil< L\le 5\left\lceil\frac{n}{N_1}\right\rceil.
\end{equation}
By \eqref{B1}, $\frac{4n}{N_1}<N_2$, which together with \eqref{UK5} gives $L\le N_2$, as required above.
We conclude from \eqref{UL0} and \eqref{B6} 
\begin{equation}
\label{UL4}
 \frac{N_2}{10}\left\lceil\frac{n}{N_1}\right\rceil^{-1}<|D_{j}| \le N_2\left\lceil\frac{n}{N_1}\right\rceil^{-1}.
\end{equation}
Let for $1\le i\le N_1$ and $1\le j\le L$
\begin{equation}
\label{UK8}
\psi_{ij}(s,t) = \left\{ \begin{array}{ll}
                N_1^{1/p}N_2^{1/u}|D_{j}|^{-1/u} & \text{ if} \enspace s=i\quad\mbox{and}\quad t\in D_{j},\\
                0  & \text{ otherwise,}
                         \end{array}
                 \right.
\end{equation}
and let $\mu_{n,N_1,N_2}^{(1)}$ be the uniform distribution on the set 
$$
\{ \alpha \psi_{ij}:\,i=1, \dots, N_1,\, j=1, \dots, L ,\, \alpha =\pm 1\} \subset B_{L_p^{N_1}(L_u^{N_2}))} .
$$
Recall that by (\ref{B6}), $LN_1>4n$, so from Lemma \ref{lem:5}(ii) and relation \eqref{UL4} we conclude
\begin{eqnarray*}
e_n^{\rm avg }(J^{N_1,N_2},\mu_{n,N_1,N_2}^{(1)},L_q^{N_1}(L_v^{N_2}))
&\ge& \frac12\big\| J^{N_1,N_2} \psi_{1,1}\big\|_{L_q^{N_1}(L_v^{N_2})}=\frac12 N_1^{1/p-1/q}N_2^{1/u-1/v}|D_1|^{1/v-1/u} 
\\
&\ge &    c N_1^{1/p-1/q}\left\lceil\frac{n}{N_1}\right\rceil^{1/u-1/v},
\end{eqnarray*}
thus \eqref{UL2}.

To prove \eqref{L5}, we set 
\begin{equation}
\label{UK6}
M=\left\lfloor\frac{4n}{N_2}\right\rfloor +1,
\end{equation}
so similarly to the above,
\begin{equation}
\label{M2}
\left\lceil\frac{4n}{N_2}\right\rceil< M\le 5\left\lceil\frac{n}{N_2}\right\rceil
\end{equation}
and $M\le N_1$. Now define for $1\le i\le M$ and $1\le j\le N_2$
\begin{equation}
\label{N1}
\psi_{ij}(s,t) = \left\{ \begin{array}{ll}
                N_1^{1/p}M^{-1/p}, & \text{ if} \enspace s=i\quad\mbox{and}\quad t=j\\
                0  & \text{ otherwise.}
                         \end{array}
                 \right.
\end{equation}
Let $(\varepsilon_{ij})_{i=1,j=1}^{M,N_2}$ be independent 
symmetric Bernoulli random variables and let $\mu_{n,N_1,N_2}^{(2)}$ be the distribution of 
$
\sum_{i=1}^M \sum_{j=1}^{N_2} \varepsilon_{ij}\psi_{ij}.
$
Then $\mu_{n,N_1,N_2}^{(2)}$ is concentrated on $B_{L_p^{N_1}(L_u^{N_2})}$.
Since by (\ref{M2}), $MN_2>4n$, we can apply Lemma \ref{lem:5}. So let $\mathcal{K}$ be any subset of 
$\{(i,j)\,:\,1\le i\le M,\, 1\le j\le N_2\}$ with $|\mathcal{K}|\ge MN_2-2n$. Then 
\begin{equation}
\label{Umc-eq:5}
|\mathcal{K}|\ge\frac12 MN_2.
\end{equation}
For $1\le i\le M$ let 
\begin{equation}
\label{M1}
\mathcal{K}_i=\{1\le j\le N_2\,:\, (i,j)\in \mathcal{K}\},\quad I := \left\{ 1\le i\le M \, :\, |\mathcal{K}_i| \geq \frac{N_2}{4} \right\} .
\end{equation}
Then $|I| \geq \frac{M}{4}$ and 
we get from \eqref{M2} and \eqref{M1}
\begin{eqnarray*}
{\mathbb E}\, \bigg\|   \sum_{(i,j)\in \mathcal{K}} \varepsilon_{ij} J^{N_1,N_2}\psi_{ij}  \bigg\|_{L_q^{N_1}(L_v^{N_2})}
 &\ge& {\mathbb E}\, \bigg\|  \sum_{i\in I} \sum_{j\in \mathcal{K}_i} \varepsilon_{ij}\psi_{ij} \bigg\|_{L_q^{N_1}(L_v^{N_2})}\ge 4^{-1/v}N_1^{1/p}M^{-1/p}|I|^{1/q}N_1^{-1/q}
\\
&\ge&  cN_1^{1/p-1/q}M^{1/q-1/p}\ge cN_1^{1/p-1/q}\left\lceil\frac{n}{N_2}\right\rceil^{1/q-1/p}
\end{eqnarray*}
and from Lemma \ref{lem:5} (i)
\begin{eqnarray*}
\lefteqn{e_n^{\rm avg }(J^{N_1,N_2},\mu_{n,N_1,N_2}^{(2)},L_q^{N_1}(L_v^{N_2}))}
\notag\\ 
&\ge& \frac12\min_{|\mathcal{K}|\ge M N_2-2n}{\mathbb E}\, \Bigg\|   \sum_{(i,j)\in \mathcal{K}} \varepsilon_{ij} J^{N_1,N_2}\psi_{ij}  \Bigg\|_{L_q^{N_1}(L_v^{N_2})}
\ge cN_1^{1/p-1/q}\left\lceil\frac{n}{N_2}\right\rceil^{1/q-1/p},
\end{eqnarray*}
thus \eqref{L5}.
 
We derive relation \eqref{N3} directly  from \eqref{UL2} and \eqref{UL1} from \eqref{L5}. Setting $n_1=\lceil c_0N_1N_2\rceil-1$ and recalling from \eqref{B1} that $1<c_0N_1N_2$, we get
$
\frac{c_0}{2}N_1N_2\le n_1< c_0N_1N_2.
$
Consequently, 
\begin{equation}
\label{M4}
\frac{c_0}{2}N_1\le \left\lceil\frac{n_1}{N_2}\right\rceil< (c_0+1)N_1, \quad   \frac{c_0}{2}N_2\le \left\lceil\frac{n_1}{N_1}\right\rceil< (c_0+1)N_2.
\end{equation}
We set 
\begin{eqnarray*}
\mu_{n,N_1,N_2}^{(3)}=\mu_{n_1,N_1,N_2}^{(1)},\quad
\mu_{n,N_1,N_2}^{(4)}=\mu_{n_1,N_1,N_2}^{(2)}.
\end{eqnarray*}
Furthermore, since $n<c_0N_1N_2$ we have $n\le n_1$, hence by monotonicity, \eqref{UL2}, and \eqref{M4}
\begin{eqnarray*}
\lefteqn{e_n^{\rm avg }(J^{N_1,N_2},\mu_{n,N_1,N_2}^{(3)},L_q^{N_1}(L_v^{N_2}))= e_n^{\rm avg }(J^{N_1,N_2},\mu_{n_1,N_1,N_2}^{(1)},L_q^{N_1}(L_v^{N_2}))}
\\
&\ge& e_{n_1}^{\rm avg }(J^{N_1,N_2},\mu_{n_1,N_1,N_2}^{(1)},L_q^{N_1}(L_v^{N_2})) 
\ge  c N_1^{1/p-1/q}\left\lceil\frac{n_1}{N_1}\right\rceil^{1/u-1/v}
\ge c  N_1^{1/p-1/q}N_2^{1/u-1/v},
\end{eqnarray*}
thus \eqref{N3}.
Similarly, from \eqref{L5} 
\begin{eqnarray*}
\lefteqn{e_n^{\rm avg }(J^{N_1,N_2},\mu_{n,N_1,N_2}^{(4)},L_q^{N_1}(L_v^{N_2}))= e_n^{\rm avg }(J^{N_1,N_2},\mu_{n_1,N_1,N_2}^{(2)},L_q^{N_1}(L_v^{N_2}))}
\\
&\ge& e_{n_1}^{\rm avg }(J^{N_1,N_2},\mu_{n_1,N_1,N_2}^{(2)},L_q^{N_1}(L_v^{N_2})) 
\ge cN_1^{1/p-1/q}\left\lceil\frac{n_1}{N_2}\right\rceil^{1/q-1/p}\ge c,
\end{eqnarray*}
which is \eqref{UL1}.

For the proof of inequalities \eqref{UL3} and \eqref{UQ3} we can assume $n\ge N_1$, because for $n<N_1$ the already shown relation \eqref{UL1} implies \eqref{UL3}, while \eqref{UL2} gives \eqref{UQ3}. We set 
\begin{equation}
\label{UV1}
L=4\left\lceil\frac{4n}{N_1}\right\rceil+1,
\end{equation}
hence by \eqref{B1}
$$
L\le \frac{16n}{N_1}+5\le \frac{21n}{N_1}\le N_2.
$$

To prove \eqref{UL3}, 
we use again the blocks 
$D_{j} \;( j=1, \dots, L)$  given by  (\ref{UWL7}) and define $\psi_j\in B_{L_u^{N_2}}$ by 
\begin{equation}
\label{UK9}
\psi_{j}(t) = \left\{ \begin{array}{ll}
                N_2^{1/u}|D_{j}|^{-1/u} & \text{ if} \enspace \quad t\in D_{j},\\
                0  & \text{ otherwise.}
                         \end{array}
                 \right.
\end{equation}
Let $\mu_1$ be the counting measure on $\{\pm \psi_j:\, 1\le j\le L\}\subset L_u^{N_2}$. Then we set $\mu_{n,N_1,N_2}^{(5)}=\mu_1^{N_1}$, the $N_1$-th power of $\mu_1$, considered as a measure on $L_p^{N_1}(L_u^{N_2})$. This measure has its support in $B_{L_p^{N_1}(L_u^{N_2})}$, and with $J^{N_2}:L_u^{N_2}\to L_v^{N_2}$ being the identical embedding,  Corollary 2.4 of \cite{Hei23a} gives 
\begin{equation}
\label{B9}
e_n^{\rm avg }(J^{N_1,N_2},\mu_{n,N_1,N_2}^{(5)},L_q^{N_1}(L_v^{N_2}))\ge 2^{-1-1/q}e_{\left\lceil\frac{4n}{N_1}\right\rceil}^{\rm avg }(J^{N_2},\mu_1,L_v^{N_2}).
\end{equation}
By Lemma \ref{lem:5}(ii) with $\bar{n}=L$, \eqref{UL0}, and \eqref{UV1}
\begin{eqnarray*}
e_{\left\lceil\frac{4n}{N_1}\right\rceil}^{\rm avg }(J^{N_2},\mu_1,L_v^{N_2})\ge \frac{1}{2}\|J^{N_2}\psi_1\|_{L_v^{N_2}}=\frac{1}{2}N_2^{1/u-1/v}|D_{j}|^{1/v-1/u}
\ge c \left\lceil\frac{n}{N_1}\right\rceil^{1/u-1/v},
\end{eqnarray*}
which together with  \eqref{B9} gives \eqref{UL3}.

Finally we turn to \eqref{UQ3}, where we set 
\begin{equation}
\label{UP3}
\psi_j=N_1^{1/p}\chi_{D_j}\in L_u^{N_2} \quad (j=1,\dots,L),
\end{equation}
with $D_{j}$ given by  (\ref{UWL7}) and $L$ by  (\ref{UV1}). Let $(\varepsilon_j)_{j=1}^L$ be independent symmetric Bernoulli random variables and let   $\mu_1$ be the distribution of  $ \sum_{j=1}^L  \varepsilon_j\psi_j$. We define a measure $\mu_{n,N_1,N_2}^{(6)}$ on $B_{L_p^{N_1}(L_u^{N_2})}$ as follows: Let $\Phi_k:L_u^{N_2}\to L_p^{N_1}(L_u^{N_2})$ be the identical embedding onto the $k$-th component of the space $L_p^{N_1}(L_u^{N_2})$, that is, for $g\in L_u^{N_2}$, $\Phi_k(g)=f$, with  $f(k,j)=g(j)$ and $f(i,j)=0$ for $i\ne k$.
We define the measure $\mu_{n,N_1,N_2}^{(6)}$ on $L_p^{N_1}(L_u^{N_2})$ by setting for a set $C\subset L_p^{N_1}(L_u^{N_2})$
\begin{equation}
\label{B8}
\mu_{n,N_1,N_2}^{(6)}(C)=N_1^{-1}\sum_{i=1}^{N_1}\mu_1(\Phi_i^{-1}(C)),
\end{equation}
thus by \eqref{UP3}, $\mu_{n,N_1,N_2}^{(6)}$ is of finite support in $B_{L_p^{N_1}(L_u^{N_2})}$. Now Corollary  2.6 of \cite{Hei23a} yields
\begin{equation}
\label{UP2}
e_n^{\rm avg-non}(J^{N_1,N_2},\mu_{n,N_1,N_2}^{(6)},L_q^{N_1}(L_v^{N_2}))\ge \frac{1}{2} N_1^{-1/q}e_{\left\lceil\frac{2n}{N_1}\right\rceil}^{\rm avg-non }(J^{N_2},\mu_1,L_v^{N_2}).
\end{equation}
By Lemma \ref{lem:5}(i), \eqref{UL0}, \eqref{UV1}, and \eqref{UP3}
\begin{eqnarray*}
\lefteqn{e_{\left\lceil\frac{2n}{N_1}\right\rceil}^{\rm avg-non}(J^{N_2},\mu_1,L_v^{N_2})\ge e_{\left\lceil\frac{2n}{N_1}\right\rceil}^{\rm avg }(J^{N_2},\mu_1,L_v^{N_2})}
\\
&\ge& \frac{1}{2}\min\bigg\{{\mathbb E}\,\Big\|\sum_{i\in I}\varepsilon_i J^{N_2}\psi_i\Big\|_{L_v^{N_2}}:\,I\subseteq\{1,\dots,L\},\,|I|\ge L-2\left\lceil\frac{2n}{N_1}\right\rceil\bigg\}
\ge  cN_1^{1/p}.
\end{eqnarray*}
Inserting this into \eqref{UP2} finally yields \eqref{UQ3}.
\end{proof}

\begin{theorem}
\label{theo:1}
Let $1 \leq p,q,u,v \le \infty $. Then there exist constants $0<c_0<1$, $c_1,\dots,c_6>0$,  such that for all $n,N_1,N_2\in{\mathbb{N}}$ with 
$n < c_0N_1N_2$ the following hold: 
If  $p\ge q$ or $u\le v$, then  
\begin{eqnarray}
\label{A2}
&&c_1N_1^{\left(1/p-1/q\right)_+}N_2^{\left(1/u-1/v\right)_+}
\le 
e_n^{\rm ran }\Big(J^{N_1,N_2},B_{L_p^{N_1}(L_u^{N_2})},L_q^{N_1}(L_v^{N_2})\Big) \nonumber\\[.2cm]
&\le& e_n^{\rm ran-non}\Big(J^{N_1,N_2},B_{L_p^{N_1}(L_u^{N_2})},L_q^{N_1}(L_v^{N_2})\Big)
\le 
c_2N_1^{\left(1/p-1/q\right)_+}N_2^{\left(1/u-1/v\right)_+}.
\end{eqnarray}
If $p<q$ and $u>v$, then  
\begin{eqnarray}
\label{A3}
&&c_3N_1^{1/p-1/q}\Bigg(\left\lceil\frac{n}{N_1}\right\rceil^{1/u-1/v}+ \left\lceil\frac{n}{N_2}\right\rceil^{1/q-1/p}\Bigg)
\le e_n^{\rm ran }\Big(J^{N_1,N_2},B_{L_p^{N_1}(L_u^{N_2})},L_q^{N_1}(L_v^{N_2})\Big) 
\nonumber\\[.2cm]
&\le& c_4N_1^{1/p-1/q}\Bigg(\left\lceil\frac{n}{N_1\log(N_1+N_2)}\right\rceil^{1/u-1/v}+ \left\lceil\frac{n}{N_2\log(N_1+N_2)}\right\rceil^{1/q-1/p}\Bigg)
\end{eqnarray}
and
\begin{eqnarray}
c_5N_1^{1/p-1/q}
  &\le& e_n^{\rm ran-non}\Big(J^{N_1,N_2},B_{L_p^{N_1}(L_u^{N_2})},L_q^{N_1}(L_v^{N_2})\Big) 
\le c_6N_1^{1/p-1/q}.
\label{A4}
\end{eqnarray}
\end{theorem}
\begin{proof}
The upper bounds in \eqref{A2} and \eqref{A4} are a consequence of \eqref{WN2}, just using the zero algorithm. If $n<6(N_1+N_2)\left\lceil c(1)\log(N_1+N_2)\right\rceil$, where $c(1)>1$ is the constant $c_1$ from Proposition \ref{pro:5}, the upper bound of \eqref{A3} follows from \eqref{WN2}, as well. Now assume
\begin{equation}
\label{C1}
n\ge 6(N_1+N_2)\left\lceil c(1)\log(N_1+N_2)\right\rceil
\end{equation}
and set 
\begin{equation}
\label{M9}
m=\left\lceil c(1)\log(N_1+N_2)\right\rceil, \quad\tilde{n}=\left\lfloor\frac{n}{6m}\right\rfloor.
\end{equation}
We use Proposition \ref{pro:5} with $\tilde{n} $ instead of $n$, so by \eqref{WM1} and \eqref{M9}
\begin{equation}
{\rm card}(A_{\tilde{n},\omega}^2)\le (m+1)\tilde{n}+mN_1+N_2\le 2m\tilde{n}+m(N_1+N_2)
\le 3m\tilde{n}\le \frac{n}{2},
\end{equation}
consequently,
$
{\rm card}(A_{\tilde{n},\omega}^3)\le 2\,{\rm card}(A_{\tilde{n},\omega}^2)\le n
$
and therefore
\begin{eqnarray}
\label{C2}
e_n^{\rm ran }\Big(J^{N_1,N_2},B_{L_p^{N_1}(L_u^{N_2})},L_q^{N_1}(L_v^{N_2})\Big) 
&\le & 
 cN_1^{1/p-1/q}\Bigg(\left\lceil\frac{\tilde{n}}{N_1}\right\rceil^{1/u-1/v}+ \left\lceil\frac{\tilde{n}}{N_2}\right\rceil^{1/q-1/p}\Bigg).
\end{eqnarray}
Furthermore, we obtain from \eqref{C1} and \eqref{M9}
\begin{eqnarray}
\left\lceil\frac{\tilde{n}}{N_i}\right\rceil 
&\ge& \frac{cn}{N_i m}
\ge \frac{cn}{N_i\log(N_1+N_2)}\ge c\left\lceil \frac{n}{N_i\log(N_1+N_2)}\right\rceil \quad (i=1,2).\label{C3}
\end{eqnarray}
Combining \eqref{C2} and \eqref{C3} proves the upper bound in \eqref{A3}.

Now we prove the lower bounds in \eqref{A2}--\eqref{A4}. We use Lemma \ref{Ulem:5} and Proposition \ref{pro:1}.  We assume $n<\frac12 c(0) N_1N_2$, where $c(0)$ stands for the constant $c_0$ from Proposition \ref{pro:1}. We start with \eqref{A2} and assume  first that $p\ge q$ and $u\le v$. Setting $n_1=\left\lceil\frac12 c(0) N_1N_2\right\rceil-1$, we have $n\le n_1$ and 
$$
 \frac14 c(0) N_1N_2\le n_1< \frac12 c(0) N_1N_2,
$$
therefore, by \eqref{UL3} and monotonicity of the $n$-th minimal errors
\begin{eqnarray*}
\lefteqn{e_n^{\rm ran }\Big(J^{N_1,N_2},B_{L_p^{N_1}(L_u^{N_2})},L_q^{N_1}(L_v^{N_2})\Big)\ge e_{n_1}^{\rm ran }\Big(J^{N_1,N_2},B_{L_p^{N_1}(L_u^{N_2})},L_q^{N_1}(L_v^{N_2})\Big)}
\\
&\ge& \frac12 e_{2n_1}^{\rm avg }(J^{N_1,N_2},\mu_{2n_1,N_1,N_2}^{(5)},L_q^{N_1}(L_v^{N_2}))
 \ge c \left\lceil\frac{2n_1}{N_1}\right\rceil^{1/u-1/v}\ge c N_2^{1/u-1/v},
\label{N6}
\end{eqnarray*}
which is the lower bound of \eqref{A2} for the case $(p\ge q)\wedge(u\le v)$. If  $(p\ge q)\wedge(u> v)$, the lower estimate of relation \eqref{A2} follows from \eqref{UL1}, while if $(p< q)\wedge(u\le v)$, it is a consequence of \eqref{N3}. Relations \eqref{UL2} and \eqref{L5} together give the lower bound in \eqref{A3}, and  \eqref{UQ3} implies the one in \eqref{A4}. 
\end{proof}
With  $c_0$  from Theorem \ref{theo:1} and $N_1=N_2=\left\lfloor c_0^{-1/2}n^{1/2}\right\rfloor+1$, we obtain
\begin{corollary}
\label{cor:1}
There are constants $c_1,c_2>0$ such that for each $n\in{\mathbb N}$ there exist $N_1,N_2\in{\mathbb N}$ such that
\begin{equation*}
c_1n^{1/2}(\log(n+1))^{-1}\le \frac{e_n^{\rm ran-non }(J^{N_1,N_2},B_{L_1^{N_1}(L_\infty^{N_2})},L_\infty^{N_1}(L_1^{N_2}))}{e_n^{\rm ran }(J^{N_1,N_2},B_{L_1^{N_1}(L_\infty^{N_2})},L_\infty^{N_1}(L_1^{N_2}))}\le c_2n^{1/2}.
\end{equation*}
\end{corollary}
It is not known though if the exponent $1/2$ is the largest possible among all linear problems. More precisely, is $\sup \Gamma > 1/2$, where
$$
\Gamma=\left\{\gamma>0:\, \exists c>0\, \forall n\in{\mathbb N} \,\exists \text{ a linear problem } \mathcal{P}_n \text{ with } 
\frac{e_n^{\rm ran-non }(\mathcal{P}_n)}{e_n^{\rm ran }(\mathcal{P}_n)}\ge cn^\gamma
\right\} \,?
$$

\end{document}